\RequirePackage{fix-cm}
\documentclass[smallextended]{svjour3}       
\smartqed  
\usepackage{graphicx}
\usepackage{a4wide}

\usepackage{tikz}
\tikzstyle{every picture}=[line cap=round,line join=round,line width=.7pt,minimum size=3pt,every label/.append style={font=\small}, label distance=2pt]
\tikzstyle{EDR}=[draw=red,line width=1pt,preaction={clip, postaction={pattern=north west lines, pattern color=red}}]
\tikzstyle{EDB}=[draw=red,line width=1pt,preaction={clip, postaction={pattern=north west lines, pattern color=blue}}]
\tikzstyle{EDG}=[draw=red,line width=1pt,preaction={clip, postaction={pattern=north west lines, pattern color=green}}]
\usepackage{caption}
\usepackage{color}

\usepackage{graphicx}
\usepackage{enumerate,color}
\usepackage{hyperref}
\usepackage{array}
\usepackage{amsfonts}
\usepackage{psfrag}
\usepackage{pgf,tikz}
\usetikzlibrary{arrows}
\usetikzlibrary{patterns}
\usepackage{mathtools}

\newcommand{\e}{\varepsilon}

\newcommand{\modulo}[1]{{\left|#1\right|}}

\newcommand{\R}{\mathbb{R}}

\newcommand{\N}{\mathbb{N}}

\renewcommand{\epsilon}{\varepsilon}
\renewcommand{\phi}{\varphi}
\renewcommand{\theta}{\vartheta}

\begin{document}

\title{Fujita modified exponent for scale invariant damped semilinear wave equations. 
}


\author{Felisia Angela Chiarello\and
        Giovanni Girardi\and Sandra Lucente
}


\institute{Felisia Angela Chiarello \at
              Department of Mathematical Sciences “G. L. Lagrange”, Politecnico di Torino, Corso Duca degli Abruzzi 24, 10129 Torino, Italy. Supported by  “Compagnia di San Paolo” (Torino, Italy) \\
              \email{felisia.chiarello@polito.it}           
           \and
           Giovanni Girardi \at
             Department of Mathematics, University of Bari, Via Orabona n. 4, Italy.\\
             \email{giovanni.girardi@uniba.it}    
              \and
             Sandra Lucente \at
             Dipartimento Interateneo di Fisica, University of Bari, Via Orabona n. 4, Italy. Partially supported by PRIN 2017- linea Sud  "Qualitative and quantitative aspects of nonlinear PDEs."\\
             \email{sandra.lucente@uniba.it}   
}

\date{Received: date / Accepted: date}

\maketitle

\begin{abstract}
The aim of this paper is to prove a blow up result of the solution for a semilinear scale invariant
damped wave equation under a suitable decay condition on radial initial data. The admissible range for the power of the nonlinear term depends both on the damping coefficient and on the pointwise decay order of the initial data.  In addition we give an upper bound estimate for the lifespan of the solution.
It depends not only on the exponent of the nonlinear term, not only on the damping coefficient but also on the size of the decay rate of the initial data. 
\end{abstract}

\section{Introduction}
\label{intro}
{In the recent years, the following Cauchy problem for the wave equation with scale invariant damping spreads a new line of research on variable coefficient type equations. More precisely, we are dealing with  
	\begin{equation}\label{eq:dampedmasswave}
	\begin{cases}
	v_{tt}(t,x)-\Delta v(t,x)+\frac{\mu}{1+t}v_t(t,x)+\frac{\nu}{(1+t)^2}v(t,x)=\modulo{v(t,x)}^p, \quad &t\geq0,\quad x\in \R^n,\\
	v(0,x)=0,  \\
	v_t(0,x)=\varepsilon g(x),  
	\end{cases}
	\end{equation}
	with $n\ge 2$, $\mu,\,\nu \in \R$,  $p>1$ and $g$ a radial smooth function.
In \cite{DLR2015}, \cite{DL2015} and \cite{P2019even} and \cite{P2019odd} some results on the global existence of a solution for \eqref{eq:dampedmasswave} with non compactly initial data appeared assuming  a suitable decay behavior for $g$. Many other results concern blow-up and global existence for this equation, see \cite{PR2019} for a summary of this problem. The main point is to find a critical exponent, fixed a suitable space of data. More precisely, a level $\bar p$ is critical  if for $p >\bar  p$  one can prove that for $\epsilon>0$ sufficiently small and for any $g$ chosen in the fixed space there exists a unique global (in time) solution of the problem, and conversely if $p\in (1, \bar p)$, for any $\epsilon>0$ there exist some $g$ in this space such that the local solution cannot be prolonged over a finite time. Coming back to \eqref{eq:dampedmasswave}, in dependence on $\mu,\,\nu$ and $n,$ a competition between two critical exponents appeared.  In some cases the Strauss exponent is dominant; it is given by the wave equation theory, it will be denoted by $p_S(d)$ and it is the positive root of the quadratic equation
$$ (d-1)p^2-(d+1)p-2=0\,.$$ For other assumptions, the equation goes to an heat equation and a Fujita-type exponent 
	$p_F(h):=1+\frac{2}{h}$ appears. In all known results,  the quantities $d>1, h> 0$ depend on $\nu,\,\mu$ and $n$. Changing the space of data, a change of critical exponent may appear. The novelty of our result consists in showing that if one takes into account the decay rate of the initial data then the Fujita type exponent depends also on such decay rate. In addition we give an upper bound estimate for the lifespan of the solution, in terms of the power of the nonlinear term, the size and the growth of the initial data. Let us recall that the lifespan of the solution is a function of $\varepsilon$ which gives the maximal existence time: $$ T(\varepsilon):=\sup\,\{ T>0  \text{ such that the local solution }\, u \,\text{ to } \eqref{eq:dampedmasswave} \text{ is defined on } [0,T)\times \R^n\,\}. $$

	Finally, we will prove the following.
	\begin{theorem}\label{th:main0}
		Let $n\geq 2$. Let  $\epsilon>0$ and $g$ a radial smooth function satisfying
		\begin{equation}\label{eq:hyp}
		g(|x|)\geq \frac{M}{(1+|x|)^{\bar k+1},} 
		\text{ with } \bar{k}>-1, 
		\end{equation}
		for some $M>0$ and for any $x \in \R^n$. Assuming in addition that
		$$
		  \bar{k}+\frac{\mu}{2}>0,\quad \frac{\mu}{2}\left(\frac{\mu}{2}-1\right)\geq \nu$$
		  and $$ 1<p<p_F\left(\bar k+\frac{\mu}{2}\right),
		$$
		then the classical solution of \eqref{eq:dampedmasswave} blows up. More precisely, the  lifespan of the solution $T(\varepsilon)>0$ is finite and satisfies
		\begin{equation}\label{eq:T}
		T(\varepsilon)\leq C \epsilon^{-\frac{2(p-1)}{4-(\mu+2\bar{k})(p-1)}},
		\end{equation}
		with $C>0$, independent of $\varepsilon$. 
	\end{theorem}

\begin{remark}
Recently Ikeda, Tanaka, Wakasa in \cite{ITW2020} consider a similar question for cubic  convolution nonlinearity and a critical decay appears.
\end{remark}
\begin{remark}
	In \cite{GL2020}  we will also consider a variant of problem \eqref{eq:dampedmasswave}, in which the nonlinearity depends on $v,t,v_t$  combined in a suitable way.
\end{remark}

\begin{remark}
	The lifespan estimate for the same equation with compactly supported data and $\nu\neq\frac{\mu}{2}\left(\frac{\mu}{2}-1\right)$ has been considered in \cite{PT18}. If $\nu\leq\frac{\mu}{2}\left(\frac{\mu}{2}-1\right) $ the lifespan estimate is different from (3) due to the compactness of the support of the  initial data. 
	\end{remark}

	The paper is organized as follows: in Section 2 we give an overview of the known results and we state an auxiliary theorem; in Section 3 we prove the main results.
	
	\section{Motivations}
	\subsection{The case $\mu=2$, $\nu=0$}
	Let us start with a quite simple case}
\begin{equation}
\label{eq:mu2}
\begin{cases}
v_{tt}-\Delta v+\frac{2}{1+t}v_t=\modulo{v}^p, \quad &t\geq0,\quad x\in \R^n,\\
v(0,x)=0,  \\
v_t(0,x)=\varepsilon g(x).
\end{cases}
\end{equation}
The global existence of small data solutions for this problem was first solved in \cite{DA2015} for a suitable range of $n$ and $p$. Some non-existence results were also established for $p<p_F(n):=1+\frac{2}{n}$. Except for the one-dimensional case a gap between this value and the admissible exponents in \cite{DA2015} appeared. 
In \cite{DLR2015} for dimension $n=2,3,$ this gap was covered with an unexpected result. Indeed, in that paper the Strauss exponent came into play. 
Afterwards, the global existence of small data solutions to \eqref{eq:mu2} has been proved for any $p>p_S(n+2)$ also in odd dimension $n\geq 5$ in \cite{DL2015} and in even dimension $n\geq 4$ in \cite{P2019even}.\\ 
Moreover, we know that the exponent $p_2(n):=\max\{p_S(n+2),p_F(n)\}$ is optimal; in fact, in \cite{DLR2015} the authors prove the blow up of solutions of \eqref{eq:mu2} for each $1<p\le p_2(n)$ in each dimension $n\in \N$.
In \cite{DL2015,DLR2015,P2019even}, the authors prove a global existence result not necessarily when the initial datum $g=g(x)$ has compact support. More precisely, let $n\ge3$, given a radial initial datum $g(x)=g(|x|)$ with $g\in C^1(\R)$,
for any $p>p_S(n+2)$ it is possible to choose $\bar k>0$ and $\epsilon_0>0$ such that \eqref{eq:mu2} admits a radial global solution $u\in C([0,\infty)\times \R^n)\cap C^2([0,\infty)\times \R^n\setminus \{0\})$ provided
\begin{equation}
\label{hyp:g-existence}
|g^{(h)}(r)|\leq \epsilon \langle r \rangle ^{-(\bar{k}+1+h)} \text{ for } h=0,\,1.
\end{equation}
In the present paper we discuss the dependence of $\bar k$ from $n$ and $p$.  
In  \eqref{hyp:g-existence}, the exponent $\bar{k}$ has to belong to a suitable interval $ [k_1(n,p),k_2(n,p)]$.
It is interesting to investigate the case of $\bar k\not \in [k_1(n,p), k_2(n,p)]$. In the sequel we will see that the bound $k_2(n,p)$ can be easily improved (see Remark \ref{remark:avoidK2}). On the contrary if $k<k_1(n,p)$ then a new result appears. The known situation is the following:
\begin{itemize}
	\item[-] $k_1(3,p)=\max\big\{\frac{3-p}{p-1}, \frac{1}{p-1}\big\}$ and $k_2(3,p)=2(p-1)$, see \cite{DLR2015}.
	\item[-]$k_1(n,p)= \max\big\{\frac{3-p}{p-1}, \frac{n-1}{2}\big\}$
	and $k_2(n,p)= \min\big\{\frac{(n+1)p}{2}-2,\frac{n^2-2n+13}{2(n-3)}\big\}$ if  $n\geq 5$ odd, see \cite{DL2015}.
	\item[-] $k_1(n,p)= \max\big\{\frac{3-p}{p-1},\frac{n-1}{2}\big\}$
	and $k_2(n,p)= \min\left\{\frac{(n+1)p}{2}-2, n-1\right\}$ if  $n\geq 4$, see \cite{P2019even}.
\end{itemize}  
We can write in a different way the previous conditions. 
Firstly we concentrate on the case $n=3$. 
For $p\in (1,2)$ we have $\bar k\ge \frac{3-p}{p-1}$ that is equivalent to 
\[
p\ge 1+\frac{2}{\bar k+1}=p_F(\bar k+1).
\]
From above we have $\bar k\le 2(p-1)$ that is 
\[
p\ge \frac{\bar k}{2}+1.
\]
The intersection of  $p=p_F(\bar k+1)$ and $p=1+\bar k/2$ is exactly in $\bar k= \frac{-1+\sqrt{17}}{2}$  and $p=p_S(5)$. 
We summarize the situation in Figure \ref{fi:mu2}. In the following graphs we denote in blue the zone of the known
global existence results, in red the zone of the known blow-up results. In this paper we want to cover the white zones.

\begin{center}
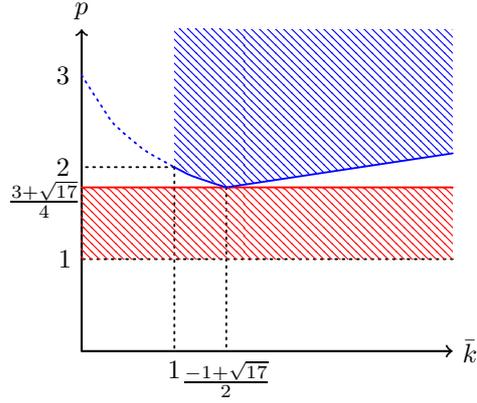

	\begin{tikzpicture}[scale=1.22]
	\draw [<->,thick] (0,3.5) node (yaxis) [above] {$p$}
	|- (4,0) node (xaxis) [right] {$\bar k$};
	\draw [dotted](4,1)-- (0,1) node [left]  (scrip) {$1$};
	\draw [dotted] (1,2)-- (0,2);
	\node[draw=none] at (-0.2, 2) {$2$};
	\node[draw=none] at (1, -0.2)   {$1$};
	\node[draw=none] at (-0.2, 3)   {$3$};
	\draw [red] (4,1.78)-- (0,1.78);
	\node [draw=none]  at (-0.4, 1.65) {$\frac{3+\sqrt{17}}{4}$};
	\draw [dotted] (1,2)-- (1,0) ;
	\draw [dotted] (1.56,1.78)-- (1.56,0) node[below]  (scrip)  {$\frac{-1+\sqrt{17}}{2}$} ;
	\draw [blue](1.56,1.78)-- (4,2.15);
	\path [EDB, draw=none] (1.56,1.78) -- (1.56,3.5) -- (4,3.5) -- (4,2.15) -- cycle;
	\draw [EDB, draw=none] (1,2) rectangle (1.56,3.5);
	\draw [EDR, draw=none] (0, 1) rectangle (4,1.78);
	\path [EDB, draw=none] (1.56,1.78) -- (1.56,2) -- (1,2)  -- cycle;
	\draw [blue] plot [ smooth] coordinates {(1,2) (1.2,1.9) (1.56, 1.78)};
	\draw [blue, dotted] plot [ smooth] coordinates {(0,3) (0.3,2.53) (0.4, 2.42) (0.5,2.33) (0.7, 2.17) (1,2) };
	\end{tikzpicture}
	\captionof{figure}{$n=3,\; \mu=2,\; \nu=0$}\label{fi:mu2}
\end{center}

\noindent
Reading \cite{DL2015} we see that the same situation appears for any odd $n\ge 5$. The critical curve
\[
p=p_F(\bar k+1)
\]
intersect the line 
$$
p=\frac{2(\bar k+2)}{n+1}
$$
in the Strauss couple 
$$ 
\left(\bar k_0,\frac{2(\bar k_0+2)}{n+1} \right)= \left(\frac{n-5+\sqrt{n^2+14n+17}}{4}, p_S(n+2) \right)\,.
$$
The only difference with the case $n=3$ is that, in the global existence zone, a bound from above appears for $p$ and this has some influence on $k_2(n,p)$. 
More precisely one can take
\[p\le \frac{n+1}{n-3},\quad \bar k\le \frac{n^2-2n+13}{2(n-3)}\,\quad \text{ if } n\ge 7,\]
and $p\le 2$, $\bar k\le 3$ if $n=5$.
Hence, the result of such paper can be represented as in Figure \ref{fi:mu2n}. Our aim is to prove blow up in the white zone below the Fujita curve.
\begin{center}
	\begin{tikzpicture}[scale=1.1]
	\draw [<->,thick] (0,3.5) node (yaxis) [above] {$p$}
	|- (7,0) node (xaxis) [right] {$\bar k$};
	\draw [dotted](7,1)-- (0,1) node [left]  (scrip) {$1$};
	\draw [dotted] (1,2)-- (0,2);
	\node[draw=none] at (-0.4, 2) {$\frac{n+5}{n+1}$};
		\node[draw=none] at (-0.4, 3.15) {$\frac{n+1}{n-3}$};
	\draw [dotted](0,3.2)-- (7.,3.15);
	\node[draw=none] at (1, -0.2)   {$\frac{n-1}{2}$};
	\draw [red] (7,1.78)-- (0,1.78);
	\node [draw=none]  at (-0.55, 1.65) {{\small $p_S(n+2)$}};
	\draw [dotted] (1,2)-- (1,0);
	\draw [dotted] (1.56,1.78)-- (1.56,0) node[below]  (scrip)  {$\bar k_0$} ;
	\draw [blue](1.56,1.78)-- (7.1,3.15);
	\path [EDB, draw=none] (1.56,1.78) -- (1.56,3.15) -- (7.1,3.15) -- (3.3,2.2) -- cycle;
	\draw [EDB, draw=none] (1,2) rectangle (1.56,3.15);
	\draw [EDR, draw=none] (0, 1) rectangle (7,1.78);
	\path [EDB, draw=none] (1.56,1.78) -- (1.56,2) -- (1,2)  -- cycle;
	\draw [blue] plot [ smooth] coordinates {(1,2) (1.2,1.9) (1.56, 1.78)};
	\draw [blue, dotted] plot [ smooth] coordinates {(0,3) (0.3,2.53) (0.4, 2.42) (0.5,2.33) (0.7, 2.17) (1,2) };
	\end{tikzpicture}
	\captionof{figure}{$n\ge 5 \text{ odd, }\; \mu=2,\; \nu=0$}\label{fi:mu2n}
\end{center}

{Even dimension is more delicate. In \cite{P2019even} the global existence result is established in the blue zone below the line $p=\frac{n+5}{n+1}$ except on the curve $p=p_F(\bar k+1)$. For convenience of the reader, we precise that in the notation of \cite{P2019even}  the role of $\bar k$ is taken by the quantity $k+\frac{n+1}{2}$.}

\begin{center}
	\begin{tikzpicture}[scale=1.22]
	\draw [<->,thick] (0,3.2) node (yaxis) [above] {$p$}
	|- (4,0) node (xaxis) [right] {$\bar k$};
	\draw [dotted](4,1)-- (0,1) node [left]  (scrip) {$1$};
	\draw [dotted] (3,2)-- (0,2);
	\node[draw=none] at (-0.4, 2) {$\frac{n+5}{n+1}$};
	\node[draw=none] at (1, -0.2)   {$\frac{n-1}{2}$};
	\node[draw=none] at ( 3, -0.25)   {$n+1$};
	\draw [dotted] (3,0)-- (3,2);
	\draw [red] (4,1.78)-- (0,1.78);
	\node [draw=none]  at (-0.6, 1.65) {{\small $p_S(n+2)$}};
	\draw [dotted] (1,2)-- (1,0);
	\draw [dotted] (1.56,1.78)-- (1.56,0) node[below]  (scrip)  {$\bar k_0$} ;
	\draw [blue](1.56,1.78)-- (4,2.15);
	\path [EDB, draw=none] (1.56,1.78) -- (1.56,2) -- (3,2)  -- cycle;
	\draw [EDR, draw=none] (0, 1) rectangle (4,1.78);
	\path [EDB, draw=none] (1.56,1.78) -- (1.56,2) -- (1,2)  -- cycle;
	\draw [blue, dotted] plot [ smooth] coordinates {(0,3) (0.3,2.53) (0.4, 2.42) (0.5,2.33) (0.7, 2.17) (1,2) (1.2,1.9) (1.56, 1.78)};
	\end{tikzpicture}
	\captionof{figure}{$n\ge 4 \text{ even, }\; \mu=2,\; \nu=0$}\label{fi:mu2ne}
\end{center}

\subsection{The case $\mu>2$ and $\nu=\frac{\mu}{2}(\frac{\mu}{2}-1)$}

In \cite{P2019even} and \cite{P2019odd} the author considers the Cauchy problem \eqref{eq:dampedmasswave} for the semilinear wave equation with scale invariant damping and mass terms, that is $\nu=\frac{\mu}{2}(\frac{\mu}{2}-1)\geq 0$. 
We see that for $\mu=2,$ it reduces to \eqref{eq:mu2}. 
Global existence of solutions to \eqref{eq:dampedmasswave} holds under the conditions
\[
\mu\in [2,M(n)], \qquad M(n)=\frac{n-1}2\left(1+\sqrt{\frac{n+7}{n-1}}\right)\,.
\]
In the even case \cite{P2019even},
the initial data satisfies \eqref{hyp:g-existence} for $\bar k\in (k_1(n,p,\mu),k_2(n,p,\mu)]$ such that
\begin{align}
&k_1(n,p,\mu) = \max\Big\{\frac{n-1}{2}, \frac{2}{p-1}-\frac{\mu}{2}\Big\}\,;\\
& k_2(n,p,\mu) = \min\Big\{n-1, \frac{n+\mu-1}{2}p-\frac{\mu+2}{2}\Big\}.\label{eq:k2conditionmueven}
\end{align}
Rewriting these conditions in term of $p$, we find that 
\[
p>p_F\left(\bar k+\frac{\mu}{2}\right), \qquad p\ge \frac{2\bar k+\mu+2}{n+\mu-1}.
\]
The intersection of the curves those define the global existence zone gives $p=p_S(n+\mu)$.
Hence, the condition $p>p_S(n+\mu)$ appears.
Moreover, another bound from above appears:
\[
p<\bar p:=\min\left\{p_F(\mu), p_F\left(\frac{n+\mu-1}{2}\right)\right\}.
\]
This means that different results for  large $\mu$ and small $\mu$ hold. This influences the positions of $k_1$ and $k_2$.
For our purpose it is sufficient to say that for $\mu\not=2$ and even $n$ the situation is similar to Figure \ref{fi:mu2ne}. More precisely, in  Figure \ref{fi:mune} $p_S(n+\mu)$ appears. The blow up result is indeed given in \cite{NPR}. The zone between $p=p_S(n+\mu)$ and $p=p_F\left(\bar k +\frac{\mu}{2}\right)$ is not covered by any known result. 

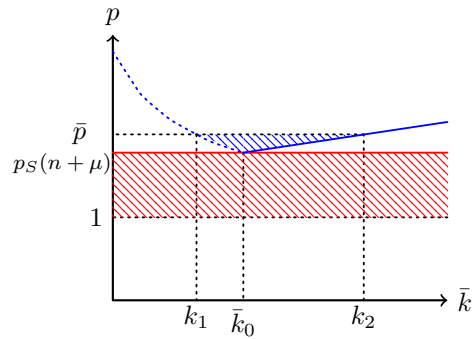
\begin{center}
	\begin{tikzpicture}[scale=1.1]
	\draw [<->,thick] (0,3.2) node (yaxis) [above] {$p$}
	|- (4,0) node (xaxis) [right] {$\bar k$};
	\draw [dotted](4,1)-- (0,1) node [left]  (scrip) {$1$};
	\draw [dotted] (3,2)-- (0,2);
	\node[draw=none] at (-0.4, 2) {$\bar p$};
	\node[draw=none] at (1, -0.2)   {$k_1$};
	\node[draw=none] at ( 3, -0.2) {$k_2$};
	\draw [dotted] (3,0)-- (3,2);
	\draw [red] (4,1.78)-- (0,1.78);
	\node [draw=none]  at (-0.6, 1.65) {{\small $p_S(n+\mu)$}};
	\draw [dotted] (1,2)-- (1,0);
	\draw [dotted] (1.56,1.78)-- (1.56,0) node[below]  (scrip)  {$\bar k_0$} ;
	\draw [blue](1.56,1.78)-- (4,2.15);
	\path [EDB, draw=none] (1.56,1.78) -- (1.56,2) -- (3,2)  -- cycle;
	\draw [EDR, draw=none] (0, 1) rectangle (4,1.78);
	\path [EDB, draw=none] (1.56,1.78) -- (1.56,2) -- (1,2)  -- cycle;
	\draw [blue, dotted] plot [ smooth] coordinates {(0,3) (0.3,2.53) (0.4, 2.42) (0.5,2.33) (0.7, 2.17) (1,2) (1.2,1.9) (1.56, 1.78)};
	\end{tikzpicture}
	\captionof{figure}{$n\ge 4 \text{ even, }\;\nu=\frac{\mu}{2}(\frac{\mu}{2}-1)\geq 0$}\label{fi:mune}
\end{center}

The corresponding global existence result for the Cauchy problem  \eqref{eq:dampedmasswave} in odd space dimension $n\geq 1$ is studied in \cite{P2019odd} for radial and small data, assuming condition \eqref{hyp:g-existence} with $\bar k \in [k_1(n,p,\mu),k_2(n,p,\mu)]$ where $k_2$ satisfies \eqref{eq:k2conditionmueven} and it holds:
\begin{align*}
&k_1(3,p,\mu)=\max\Big\{1, \frac{2}{p-1}-\frac{\mu}{2}, \frac{1}{p-1}\Big\}\,;\\
&k_1(n,p,\mu)=\max\Big\{\frac{n-1}{2}, \frac{2}{p-1}-\frac{\mu}{2}\Big\}\,, \quad n\geq 5\, \quad \mu\in [2,n-1];\\
&k_1(n,p,\mu)=\max\Big\{\frac{n-1}{2}, \frac{2}{p-1}-\frac{\mu}{2}, \frac{1}{p-1}, \Big\}\,,\, n\geq 5\quad \mu\in (n-1,M(n)].
\end{align*}
In any case the condition $p>p_F(\bar k+\frac{\mu}{2})$ appears. Hence in odd space dimension $n\ge 5$ the situation is not different from Figure \ref{fi:mune}.

Reading Theorem \ref{th:main0} in the case $\nu=\frac{\mu}{2}(\frac{\mu}{2}-1)$, it is clear that the aim of this paper is to find blowing-up solutions to \eqref{eq:dampedmasswave} even for $p>p_S(n+\mu)$ by considering initial data with slow  decay.  More precisely, let us consider
\begin{equation}\label{eq:decayl}
g(x)\simeq \frac{M}{(1+\modulo{x})^{ \bar{k}+1}}\,,\quad  \text{ for } \quad \frac{n-1}{2}<\bar k<\bar k_0\,,
\end{equation}
where $\bar k_0$ is such that 
\[
p_F\left(\bar k_0+\frac{\mu}{2}\right)=p_S(n+\mu).
\]
We will prove the blow up result in the left white side zones in Figures 1, 2, 3, 4 where  
$\bar k< \bar k_0$, $p>p_S(n+\mu)$ and $p<p_F(\bar k_0+\frac{\mu}{2})$. Under the same assumption on $g$, the quoted results assure that for $p\ge p_F(\bar k+\frac{\mu}2)$ and $p>p_S(n+\mu)$ there is global existence. Hence, $p=p_F(\bar k+\frac{\mu}{2})$ is a critical curve for the Cauchy problem \eqref{eq:dampedmasswave}, provided $\nu=\frac{\mu}{2}(\frac{\mu}{2}-1)\geq 0$.

\begin{remark}
	\label{remark:avoidK2}
	Still fixing $\nu=\frac{\nu}{2}(\frac{\nu}{2}-1)\geq 0$, let us consider $\bar{k}>\bar{k}_0$ and  $p>p_S(n+\mu)$. As discussed, the global existence results in the previous literature require $p$ above a line which depends on $\bar k$, because of a restriction of type $\bar{k}\leq k_2(n,p,\mu)$ which everytime appears.
	Actually, this restriction can be avoided; indeed, if
	the initial datum satisfies \eqref{eq:decayl} with $\bar k> k_2(n,p,\mu)$, then we can say that the
	initial datum also satisfies \eqref{hyp:g-existence} with $\bar k= k_2(n,p,\mu)$. Hence, the global existence of a solution to \eqref{eq:dampedmasswave} follows from the known results. 
\end{remark}

\begin{remark}
	For $\nu=\frac{\mu}{2}(\frac{\mu}{2}-1)\geq 0,$ Theorem \ref{th:main0} provides some new information about the solution of \eqref{eq:dampedmasswave} also when $p$ belongs to the red zone of Figure 1, 2, 3, 4, 5. In fact, for 
	\[
	p<\min\left\{p_S(n+\mu), p_F\left(\bar k+\frac{\mu}{2}\right)\right\}
	\]
	by the previous literature we know that the solution blows up in finite time, whereas Theorem \ref{th:main0} gives a life-span estimate in the case of radial initial data with non compact support, relating this estimate with the decay rate of the data.
\end{remark}

\subsection{The case $\mu=0$ and $\nu=0$}

In Figure \ref{fi:mu0} we summarize the wave equation case $\mu=\nu=0$. The red blow-up zone was covered by many authors, see \cite{Sid83} and the reference therein for the whole list of blow up results. For $\mu=\nu=0$ the global existence result has been completely solved in \cite{GLS}, where the interested reader can find a long bibliography of previous contributes. In particular the blue zone, for radial solution without compact support assumption for the initial data has been exploited by Kubo, see for example \cite{K97} and \cite{K98}. Before these papers, Takamura obtained a blow up result in the green zone. In \cite{takamura1995} the point is to find a critical decay level $k_0=\frac{2}{p-1},$ equivalently $p\leq 1+\frac{2}{k_0}.$ We underline that this a Fujita-type exponent.

\begin{center}
	\begin{tikzpicture}[scale=1.1]
	\draw [<->,thick] (0,3.5) node (yaxis) [above] {$p$}
	|- (4,0) node (xaxis) [right] {$\bar k$};
	\draw [dotted](4,1)-- (0,1) node [left]  (scrip) {$1$};
	\draw [dotted] (1,2)-- (0,2);
	\node[draw=none] at (-0.2, 2) {$3$};
	\node[draw=none] at (1, -0.15)   {$1$};
	\draw [red] (4,1.78)-- (0,1.78);
	\node [draw=none]  at (-0.35, 1.70) {$p_S(n)$};
	\draw [dotted] (1,2)-- (1,0) ;
	\draw [dotted] (1.56,1.78)-- (1.56,0) node[below]  (scrip)  {$\frac{2}{p_S(n)-1}$} ;
	\draw [EDB, draw=none] (1.56, 1.78) rectangle (4,3.5);
	\draw [EDB, draw=none] (1,2) rectangle (1.56,3.5);
	\draw [EDR, draw=none] (0, 1) rectangle (4,1.78);
	\path [EDB, draw=none] (1.56,1.78) -- (1.56,2) -- (0.9,2)-- cycle;
	\draw [blue] plot [smooth] coordinates {(0.1,2.8) (0.2,2.66) (0.3,2.53) (0.4, 2.42) (0.5,2.33) (0.6, 2.25) (0.7, 2.17) (0.8,2.1) (1,2) (1.2,1.9) (1.56, 1.78)};
	\path [EDB, draw=none] (0.1,3.5) -- (0.38,2.5)-- (0.4, 2.42)  -- (0.5,2.33) -- (0.7, 2.17) -- (1,2)--(1,3.5)-- cycle;
	\path [EDG, draw=none] (0.1, 2.8)--(0.1, 1.8)--(1.56,1.8)--(0.4,2.3)--(0.6,2.2)--(1,1.9)-- cycle;
	\end{tikzpicture}
	\captionof{figure}{$\mu=\nu=0$}\label{fi:mu0}
\end{center}
In Theorem \ref{th:main0}, we generalize Takamura's result when $\mu\neq 0$ and $\nu\leq \frac{\mu}{2}(\frac{\mu}{2}-1)$. To this aim, it is sufficient to consider a peculiar wave equation with nonlinear term having a decaying time-dependent variable coefficient. This means that we will deduce Theorem \ref{th:main0} from the following result.

\begin{theorem}\label{th:main}

	Let $n\geq 2$.  Given a smooth function $g=g(|x|)$ with $x\in \R^n$, we set $r=|x|$ and we consider $g=g(r)$ satisfying
	%
	\begin{equation}\label{eq:hyp1}
	g(r)\geq \frac{M}{(1+r)^{\bar k+1}}, 
	\text{ with } \bar k>-1,
	\end{equation}
	for some $M>0$.
	Let $u=u(t,r)$ be the radial local solution to	
	\begin{equation}
	\label{eq:radialdampedwave}
	\begin{cases}
	u_{tt}-u_{rr}-\frac{n-1}{r}u_r=(1+t)^{-\frac{\mu}{2}(p-1)}\modulo{u}^p, \quad &r>0,\\
	u(0,r)=0,  \\
	u_t(0,r)=\varepsilon g(r).
	\end{cases}
	\end{equation}
	with $p>1$ and  $p<p_{F}(\frac{\mu}{2}-1)$ if $\mu>2$. 
	Assume in addition that 
	\begin{equation}
	\label{eq:kcondition}
	-1<\bar k < \frac{2}{p-1}-\frac{\mu}{2}.
	\end{equation}
	Then, given $\epsilon>0$, the lifespan $T(\varepsilon)>0$ of classical solutions to \eqref{eq:radialdampedwave} satisfies
	\begin{equation}\label{eq:T}
	T(\varepsilon)\leq C \epsilon^{-\frac{2(p-1)}{4-(\mu+2\bar{k})(p-1)}},
	\end{equation}
	with $C>0$, independent of $\varepsilon$. 
\end{theorem}

\begin{remark}
	The assumption $p<p_{F}(\frac{\mu}{2}-1)$ if $\mu>2$ guarantees that the range of admissible $\bar k$ in \eqref{eq:kcondition} is not empty. 
\end{remark}

In the case $\mu=0$ Theorem \ref{th:main} coincides with Takamura's result in \cite{takamura1995}. In the proof of Theorem \ref{th:main}, we will follows the same approach of that paper.
\section{Proof of the main results}

\subsection{Proof of Theorem \ref{th:main}}
 We recall the crucial lemma of \cite{takamura1995}. 
\begin{lemma}\label{lemma:Tdeltam}
	Let $n\geq 2$ and $m=[n/2]$. Given a smooth function $g=g(|x|)$ with $x\in \R^n$, we set $r=|x|$ and we consider $g=g(r)$. Let us
	denote by $u^{0}(t,r)$ the solution of the free wave problem 
	\begin{align*}
	\begin{cases}
	\square u^0=0  \quad & (t,r)\in [0,\infty)\times[0,\infty)\\
	u^0(0,r)=0,~u^0_t(0,r)=g(r)\,. 
	\end{cases}
	\end{align*}
	Let $u=u(t,r)$ be a solution to 
	\begin{equation}
	\label{eq:radialF}
	u_{tt}-u_{rr}-\frac{n-1}{r}u_r=F(t,u)  
	\end{equation}
	with the initial condition 
	\begin{equation}
	\label{eq:radialFg}
	u(0,r)=0, ~ u_t(0,r)=\varepsilon g(r), \quad r\in[0,\infty).
	\end{equation}
	If $F$ is nonnegative, there exists a constant $\delta_m> 0$ such that 
	\begin{align}\label{eq:um}
	u(t,r)&\geq \varepsilon u^0(t,r)+\frac{1}{8r^m}\int_0^t d\tau \int_{r-t+\tau}^{r+t+\tau}\lambda^m F(t,u(t,\lambda))d \lambda,\\
	\label{eq:um0}
	u^0(t,r)&\geq \frac{1}{8r^m} \int_{r-t}^{r+t} \lambda^m g(\lambda) d\lambda,
	\end{align}
	provided
	\begin{equation*}
	r-t\geq \frac{2}{\delta_m}t>0.
	\end{equation*}
\end{lemma}

The constant $\delta_m$ in the previous lemma is described in  \cite[Lemma 2.5]{takamura1995}; it depends on the space dimension, 
in particular it changes accordingly with the different representations of the free wave solution in odd and even dimension.

We are ready to prove that if \eqref{eq:hyp1} holds, then the solution of \eqref{eq:radialdampedwave} blows up in finite time even for small $\varepsilon.$ 

Let us fix $\delta>0$; we define a blow-up set,
\begin{equation}
\label{eq:sigma}
\Sigma_\delta =\Big\{(t,r)\in(0,\infty)^2: r-t\geq \max\left\{\frac{2}{\delta_m}t,\delta\right\}\Big\},
\end{equation}
where $\delta_m>0$ is the constant given in Lemma \ref{lemma:Tdeltam}. Combining the assumption \eqref{eq:radialdampedwave} with the formulas \eqref{eq:um} and \eqref{eq:um0}, 
for any $(t,r)\in\Sigma_\delta,$ it holds
\begin{equation*}
u(t,r)\geq \e u^0(t,r)\geq \frac{\e}{8r^m}\int_{r-t}^{r+t}\lambda^m g(\lambda)d\lambda\geq \frac{M\e}{8r^m}\int_{r-t}^{r+t}\lambda^m (1+\lambda)^{-(\bar k +1)}d\lambda\,.
\end{equation*}
Then, \eqref{eq:sigma} implies that 
\begin{align*}
u(t,r)&\geq 
\frac{M\e}{8r^m}\left(\frac{1+\delta}{\delta}\right)^{-(\bar k+1)}\int_{r-t}^{r+t}\lambda^{m-(\bar k+1)} d\lambda\\
&\geq \frac{M\e}{8r^m}\left( \frac{1+\delta}{\delta}\right)^{-(\bar k+1)}(r+t)^{-(\bar k+1)}\int_{r-t}^{r+t}\lambda^{m} d\lambda\geq \frac{M\e}{8}\left( \frac{1+\delta}{\delta}\right)^{-(\bar k+1)}
\frac{(r-t)^{m}2t}{r^m(r+t)^{\bar k+1}}\,.
\end{align*}
Since $(t,r)\in \Sigma_\delta$, we have 
\begin{equation*}
u(t,r)\geq \frac{C_0t^{m+1}}{r^m(r+t)^{\bar k+1}}\,,
\end{equation*}
where we set 
\begin{equation}\label{eq:Czero}
C_0=\e\frac{2^{m-2}M}{\delta_m^{m}}\left(\frac{\delta}{1+\delta}\right)^{\bar k+1}>0.
\end{equation}
Now we assume an estimate of the form 
\begin{equation}
\label{eq:step0}
u(t,r)\geq \frac{Ct^a}{r^m(r+t)^b} \text{ for } (t,r)\in \Sigma_\delta,
\end{equation}
where $a$, $b$, and $C$ are positive constant. In particular, \eqref{eq:step0} holds true for $a=m+1$, $b=\bar k+1$ and $C=C_0$.

\noindent
Being $g\ge 0$, from \eqref{eq:um0} we deduce $u^0\ge 0$. Combining \eqref{eq:um} and \eqref{eq:step0}, for $(t,r)\in \Sigma_\delta$, we get
\begin{align}
\label{eq:mainstepproof}u(t,r)&\geq \frac{1}{8r^m}\int_0^t d\tau \int_{r-t+\tau}^{r+t-\tau}\frac{\lambda^m}{(1+\tau)^{\frac{\mu}{2}(p-1)}}|u(\tau, \lambda)|^p d\lambda\\
\notag&\geq \frac{C^p}{8r^m}\int_0^t \frac{\tau^{pa}}{(1+\tau)^{\frac{\mu}{2}(p-1)}}d\tau \int_{r-t+\tau}^{r+t-\tau}\lambda^{m(1-p)}(\lambda+\tau)^{-pb} d\lambda\\
& \notag\geq\frac{C^p}{8r^m(r+t)^{pb+m(p-1)}}\int_0^t \frac{\tau^{pa}}{(1+\tau)^{\frac{\mu}{2}(p-1)}}d\tau \int_{r-t+\tau}^{r+t-\tau}d\lambda\\
& \notag\geq\frac{C^p}{4r^m(r+t)^{pb+m(p-1)}}\int_0^t \frac{(t-\tau)}{(1+\tau)^{\frac{\mu}{2}(p-1)}}\tau^{pa}d\tau.
\end{align}
By means of integration by parts, we obtain 
\begin{equation*}
\int_0^t \frac{(t-\tau)\tau^{pa}}{(1+\tau)^{\frac{\mu}{2}(p-1)}}d\tau\geq 
\frac{1}{(1+t)^{\frac{\mu}{2}(p-1)}}\int_0^t (t-\tau)\tau^{pa}d\tau
\geq 
\frac{1}{(1+t)^{\frac{\mu}{2}(p-1)}}
\frac{t^{pa+2}}{(pa+1)(pa+2)}.
\end{equation*}
While searching a finite lifespan of a solution, it is not restrictive to assume $t>1$. We have
\begin{equation}\label{eq:step1}
\int_0^t \frac{(t-\tau)\tau^{pa}}{(1+\tau)^{p-1}}d\tau\geq \frac{t^{p(a-\frac{\mu}{2})+2+\frac{\mu}{2}}}{2^{p-1}(pa+1)(pa+2)}.
\end{equation}
Let $(t,r)\in \Sigma_\delta$, from \eqref{eq:step0}-\eqref{eq:step1}, we can conclude
\begin{equation}
\label{eq:generalestimate}
u(t,r)\geq \frac{C^*t^{a^*}}{r^m(r+t)^{b^*}} \text{ for }(t,r)\in\Sigma_\delta,
\end{equation}
with 
\begin{equation*}
a^*=p\Big(a-\frac{\mu}{2}\Big)+2+\frac{\mu}{2}, \hspace{3em} b^*=pb+m(p-1), \hspace{3em} C^*=\frac{(C/2)^p}{2(pa+2)^2}.
\end{equation*}
Let us define the sequences $\{a_k\},\,\{b_k\},\,\{C_k\}$ for $k\in \N$ by 
\begin{align}
&a_{k+1}=p\Big(a_k-\frac{\mu}{2}\Big)+2+\frac{\mu}{2}, \hspace{2em} a_1=m+1,\\
&b_{k+1}=pb_k+m(p-1),  \hspace{1em} b_1=\bar k+1,\\
&C_{k+1}=\frac{(C_k/2)^p}{2(pa_k+2)^2}, \hspace{1em} C_1=C_0,\label{eq:sequences}
\end{align}
where $C_0$ is defined by \eqref{eq:Czero}. Hence, we have
\begin{align}
&a_{k+1}=p^k\left(m+1-\frac{\mu}{2}+\frac{2}{p-1}\right)+\frac{\mu}{2}-\frac{2}{p-1},\label{eq:sequ1}\\
&b_{k+1}=p^k(\bar k+1+m)-m,\label{eq:sequ2} \\
& \label{eq:Ck+1first} C_{k+1}\geq K\frac{C_k^p}{p^{2k}}
\end{align}
for some constant $K=K(p,\mu, m)>0$ independent of $k$.
The relation \eqref{eq:Ck+1first} implies  that for any $k\geq 1$ it holds
\begin{align}
&C_{k+1}\geq \exp\left(p^{k}\left(\log(C_0)-S_p(k)\right)\right),
\label{eq:sequ3}\\
&S_p(k)=\Sigma_{j=0}^k d_j,\\& d_0=0 \text{ and }  d_j=\frac{j\log(p^2)-\log K}{p^{j}} \text{ for } j\geq 1.
\end{align}
We note that $d_j>0$ for sufficiently large $j.$ Since $\lim_{j\to\infty}d_{j+1}/d_j=1/p,$ the sequence $ S_p(k)$ converges for $p>1$ by using the ratio criterion for series with positive terms. Hence, there is a positive constant $S_{p,K}\ge S_p(k)$ for any $k\in \N$, so that
\begin{equation}
\label{eq:Ck+1}
C_{k+1}\geq \exp(p^{k}(\log (C_0)-S_{p,K})).
\end{equation}
Therefore, by \eqref{eq:generalestimate}, \eqref{eq:sequ1}- \eqref{eq:sequ3}, we obtain 
\begin{equation}
\label{eq:finaluestimate}
u(r,t)\geq \frac{(r+t)^m}{r^m t^{-\frac{\mu}{2}+\frac{2}{p-1}}}\exp(p^k J(t,r)),
\end{equation}
where 
\[
J(t,r):=\log(C_0)-S_{p,K}+\Big(m+1-\frac{\mu}{2}+\frac{2}{p-1}\Big)\log t-(\bar k+1+m)\log(r+t).
\]
Thus if we prove that there exists $(t_0,r_0)\in \Sigma_\delta $ such that $J(t_0,r_0)>0$, then we can conclude that the solution to \eqref{eq:radialdampedwave} blows up in finite time, in fact
\[ 
u(t_0,r_0)\to \infty  \text{ for } k\to \infty. 
\]
By the definition of $J=J(t,r)$, we find that $J(t,r)>0$ if 
\begin{align*}
\Big(\frac{2}{p-1}-\frac{\mu}{2}-\bar k\Big)\log t &> \log\Big(\frac{e^{S_{p,K}}}{C_0}\Big(2+\frac{r-t}{t}\Big)^{\bar k+1+m}\Big).
\end{align*}
In particular, we can take $(t, r)=(t, t+\max\{\frac{2t}{\delta_m},\delta\})\in \Sigma_\delta$; then, it is enough to prove that
\begin{align*}
\Big(\frac{2}{p-1}-\frac{\mu}{2}-\bar k\Big)\log t >\log\Big(\frac{e^{S_{p,K}}}{C_0}\Big(2+\frac{2}{\delta_m}\Big)^{\bar k+1+m}\Big).
\end{align*}

Now, the crucial assumption \eqref{eq:hyp1} comes into play. The coefficient in the left side is positive and by using 
\eqref{eq:Czero} we find that $J(t,r)>0$ provided 
\begin{equation}\label{eq:life}
t> C \epsilon^{-\left(\frac{2}{p-1}-\frac{\mu}{2}-\bar k\right)^{-1}},
\end{equation}
where 
\[
C=\Big(\frac{e^{S_{p,K}}\delta_m^{m}}{2^{m-2}M}\Big(\frac{1+\delta}{\delta}\Big)^{\bar k+1}\Big(2+\frac{2}{\delta_m}\Big)^{1+\bar k+m}\Big)^{\frac{1}{\frac{2}{p-1}-\frac{\mu}{2}-\bar k}},
\]
which is positive. As by-product, the inequality \eqref{eq:life} gives the lifespan estimate
\eqref{eq:T} and conclude the proof of Theorem 2.

\subsection{Proof of Theorem 1}
We start rewriting the Cauchy problem \eqref{eq:dampedmasswave} as a nonlinear wave equation with a time dependent potential.
Let $v=v(t,x)$ be a solution of \eqref{eq:dampedmasswave}; we define
\[
u(t,x):=(1+t)^{\frac{\mu}{2}}v(t,x).
\]
Then the function $u=u(t,x)$ is a solution of the Cauchy problem
\begin{equation}
\label{eq:nodampedwave}
\begin{cases}
u_{tt}-\Delta u=(1+t)^{-\frac{\mu}{2}(p-1)}\modulo{u}^p+(\frac{\mu}{2}(\frac{\mu}{2}-1)-\nu)\frac{u}{(1+t)^2}, \quad &t\geq0,\quad x\in \R^n,\\
u(0,x)=0, \\
u_t(0,x)=\varepsilon g(x).
\end{cases}
\end{equation}
If $g$ is radial, then $u$ is radial and it satisfies equations \eqref{eq:radialF} and \eqref{eq:radialFg} with $$F(t,u)=(1+t)^{-\frac{\mu}{2}(p-1)}\modulo{u}^p+\left(\frac{\mu}{2}\left(\frac{\mu}{2}-1\right)-\nu\right)\frac{u}{(1+t)^2}.$$

Let us fix $\delta>0$ and use the same notation of the proof od Theorem \ref{th:main}.
Since we are assuming $\frac{\mu}{2}(\frac{\mu}{2}-1)-\nu\geq 0$, by comparison lemma, see  \cite[Lemma 2.9]{takamura1995}, we deduce $u>0$ in $\Sigma_\delta.$
Then it holds 
$$  F(t,u)\geq (1+t)^{-\frac{\mu}{2}(p-1)}\modulo{u}^p;$$
hence, by formula \eqref{eq:um} in Lemma \ref{lemma:Tdeltam} we still derive the estimate \eqref{eq:mainstepproof}. Thus, the proof of Theorem \ref{th:main} guarantees the result of Theorem \ref{th:main0}.






\end{document}